\documentclass{amsart}
\usepackage{graphicx}
\usepackage{mathptmx}      

\usepackage{natbib}
\usepackage{geometry}
\usepackage{tikz}
\usepackage{url}
\usepackage{amsfonts}

\def\red#1 {\textcolor{red}{#1}}
\def\b#1 {\textcolor{blue}{#1}}
\newcommand{\A}{\mathcal A}
\newcommand{\D}{\mathcal D}

\newcommand{\Z}{\mathbb Z}
\newcommand{\s}{\sigma}

\begin{document}

\title{An algebraic view of bacterial genome evolution}
\thanks{This research was funded by Australian Research Council Future Fellowship FT100100898.}

\author{Andrew R. Francis}
\address{				Centre for Research in Mathematics,
              University of Western Sydney, Australia}
              \email{a.francis@uws.edu.au}           

\date{\today }

\maketitle

\begin{abstract}
Rearrangements of bacterial chromosomes can be studied mathematically at several levels, most prominently at a local, or sequence level, as well as at a topological level.  The biological changes involved locally are inversions, deletions, and transpositions, while topologically they are knotting and catenation.  These two modelling approaches share some surprising algebraic features related to braid groups and Coxeter groups.  
The structural approach that is at the core of algebra has long found applications in sciences such as physics and analytical chemistry, but only in a small number of ways so far in biology. And yet there are examples where an algebraic viewpoint may capture a deeper structure behind biological phenomena.  This article discusses a family of biological problems in bacterial genome evolution for which this may be the case, and raises the prospect that the tools developed by algebraists over the last century might provide insight to this area of evolutionary biology. .
\end{abstract}

\section{Introduction}

Evolutionary processes on the bacterial genome are dynamic and complex, with a tremendous range of mutation events occurring at a number of different physical scales.  Aside from point mutations at the level of nucleotides, a wide variety of {evolutionary} mechanisms involve cutting and rejoining of genetic material.
In this paper we step back from the relatively frequent rearrangements that occur at the single nucleotide level to look at these larger scale changes. We will reprise some of the mathematical approaches to their study, and show how it may be possible to view them in a single algebraic modelling framework.  

Examples of larger scale changes include \emph{deletion}, \emph{translocation}, \emph{duplication} and \emph{inversion}. These processes respectively delete a segment of DNA from the genome, relocate a segment to another region on the genome, make a copy of a segment and insert the duplicate into the genome, or invert a segment --- excise it and reinsert it with the opposite orientation.  These mutations are facilitated by the actions of enzymes that reside within the bacterial cell and are encoded by genes on the chromosome.  Of all the changes that may occur on the genome of a single celled bacterium, many may be fatal to the organism, either by disrupting some function essential to life, or by disrupting the replication process, and more generally may not be observed because the change has conferred a significant cost to their fitness.  
In that context it is remarkable that we know as much as we do.  We know, for instance, that of the mutation events listed above, inversion is most common, at least in bacteria~\citep{eisen2000evidence}.

Stepping back {from the focus on changes at the sequence level, } 
biological processes giving rise to \emph{knotting} in DNA have been observed for some time, at least as far back as 1981~\citep{liu1981knotted}. 
However their importance in chemistry had been recognised decades earlier by Frisch and Wassermann, who defined \emph{topological isomerism} in which chemically identical polymers differ only by their topology~\citep{frisch1961chemical} {(here we use the word ``topology" in the mathematical sense, so that differing by topology means differing \emph{as knots}). } 
For instance, polymer chains were able to be generated in the laboratory that formed links, or \emph{catenanes}, as they are known in the biochemical literature.  From a topological viewpoint,  the structures observed at this stage were simple Hopf links or trefoil knots, however a wider variety of knots were later found to arise in DNA (e.g.~\citet{cozzarelli1984topological,arsuaga2002knotting}).   

We will refer to the study of changes at a sequence level (such as inversion) as \emph{local} evolution, and those at a topological level as \emph{topological} evolution. 
While the local view of the bacterial genome focuses on the sequence and ignores the topology (or knotting), the topological view does exactly the opposite.  It is important to acknowledge that both viewpoints have sound justifications.  The sequence (local) view is natural because genes, or functional segments of DNA, are sequences, and they are transcribed in a linear way that does not take into account any twisting or topological characteristic of the location of the gene.  Other features encouraging a {sequence} view include the observation that genes that form part of the same metabolic pathway are often clustered in the same region of the genome~\citep{demerec1964clustering,stahl1966evolution,ballouz2010conditions}.  

On the other hand the topological view allows us to address a range of questions to do with the origins and maintenance of knotting, linking and supercoiling. 
For instance, the distribution of knots in the wild is not random, suggesting either that selection favours certain knot forms, or that the mechanism that gives rise to knotting leads to some knots more often than others.  Studying the distribution of knots has helped uncover information about action of the {site-specific recombinase } that is responsible for knotting in some circumstances (for example~\citet{Ernst-TanglesDNA-MPCPS1990}. 

\begin{figure}[ht]
\begin{center}
\includegraphics[width=84mm]{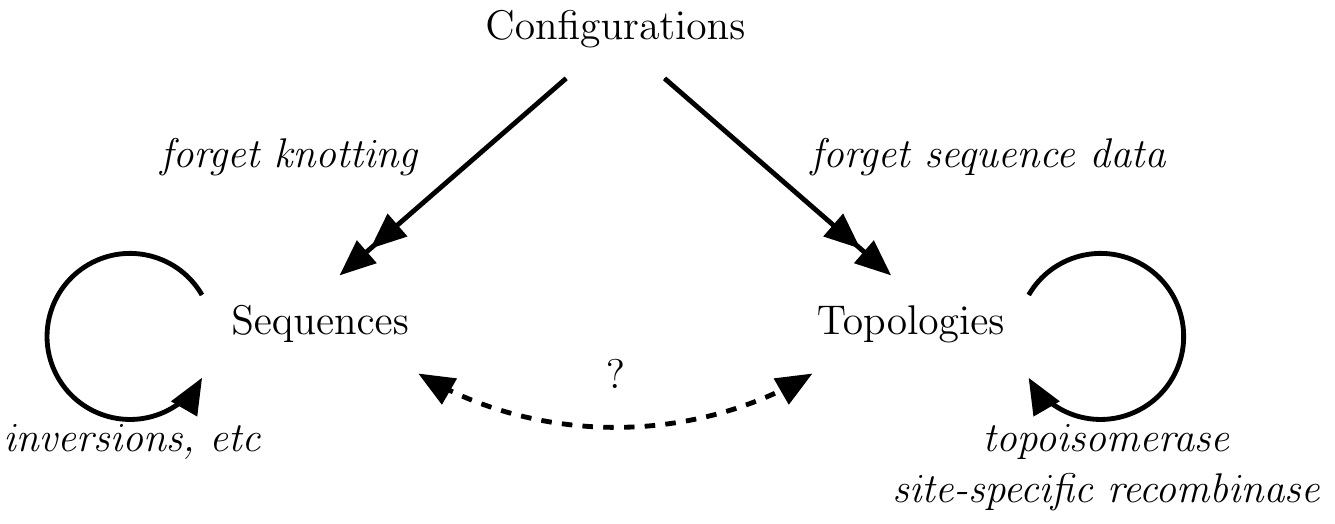}
\end{center}
\caption{Local and topological blindspots. 
}\label{fig:config}
\end{figure}

Ultimately, what we are studying in local and topological processes are two alternative projections from the actual configuration of a bacterial genome, including both nucleotide sequence data (a one-dimensional projection) and topological data (three dimensions), as in Figure~\ref{fig:config}.

Motivating both these approaches to bacterial evolution are questions about the construction of phylogenies: understanding the processes that drive the changes in the structure at a local or topological level gives information about the relationships among taxonomic units.
In this paper we will describe some algebraic structures that may provide a link between these processes that draws out the biological commonalities.  This is fertile ground for future developments.

The application of algebraic methods to biology extends beyond the use of knot theory in DNA described in this paper.  For instance, a significant body of work now applies ideas from algebraic geometry to phylogeny through the use of algebraic varieties~\citep{allman2007molecular,eickmeyer2008optimality}.  The use of varieties in biology is often termed ``algebraic biology"~\citep{pachter2005algebraic} and is closely related to the field of ``algebraic statistics"~\citep{pistone2001algebraic,drton2009lectures}.  The geometric viewpoint also has applications in viral capsid assembly~\citep{sitharam2006modeling} and RNA folding~\citep{heitsch2003rna,apostolico-finding-3D-motifs-2009,laing2011computational}, areas in which combinatorics and graph theory play a significant role, in addition to the study of radiation-induced chromosomal aberrations~\citep{sachs2002using}.  Finally, a strain of research applies group theory to problems in evolutionary biology, for instance \citet{sumner2008markov} and \citet{moulton2011butterfly} as well as the recent study of inversion distance closely related to this survey~\citep{egrinagy2013group}.

\section{Underlying biological mechanisms}

The biological mechanisms that give rise to changes such as inversion and knotting are actions of enzymes that involve cutting and rejoining DNA double-helices.  The two main families of enzymes involved in these processes are the topisomerases and the site-specific recombinases (see~\citet{yang2010topoisomerases}, for instance, for a review of these).

Topoisomerases are essential in cell replication because they cut the phosphodiester bonds between DNA base pairs, enabling untangling of the coiled or knotted DNA (see e.g~\citet{zechiedrich1995roles,hardy2004disentangling,lopez2011topo}). 
In general these do not require a specific site for cleavage, and cut at most two DNA strands (one double-helix).   Cutting one strand can allow the relaxing of super-coiling, while cutting two strands can allow the passage of one double-helix through another and aid in reducing knotting.

\begin{figure}[ht]
\begin{center}
\includegraphics[width=84mm]{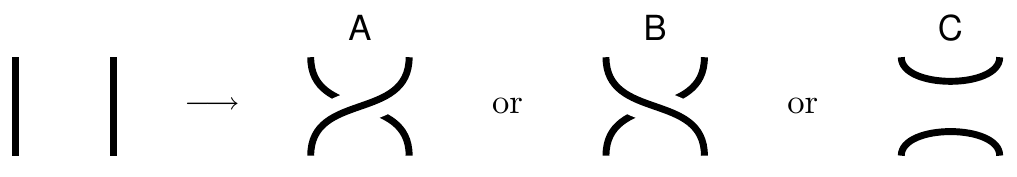}
\end{center}
\caption{Three possible transformations of the configuration of a pair of double-helical strands via cutting and rejoining double-helical strands.  Achieving any of (A), (B) or (C) from the two parallel strands shown on the left requires cutting and rejoining two double-helical strands, however configurations (A) and (B) can be reached from each other as a result of a single cut and rejoin of just one strand.  The lines in this Figure and in similar ones later in the paper represent the axis of the DNA double helix. 
}\label{fig:transf}\label{fig:dna-cutting-rejoining}
\end{figure}

In contrast, site-specific recombinases do require a specific sequence for cleavage, and act by cutting and rejoining two double-helices.
These enzymes act in a two step process, first forming a synaptic complex that involves a certain configuration of the substrate DNA, and then causing strand exchange.  
There is a rich variety of such recombinases, broadly falling into two families, the resolvase (examples include Tn3) or invertase (Gin) family, and the integrase family (phage $\lambda$, Cre, Flp, and Xer system).  While the resolvase/invertase family generally requires topological alignment of sites that are directly repeated or inverted (respectively), the integrase family are more versatile and can act on a wider range of substrate arrangements (see~\citet{crisona1999topological}, or the Introduction to~\citet{craig2002mobile} for surveys of these recombinases).  These types of action can take place as a result of cuts to a single double-helical strand at the site of recombination, or to two double-helical strands, as shown in Figure~\ref{fig:dna-cutting-rejoining}. %

\begin{figure}[ht]
\includegraphics[width=8.4cm]{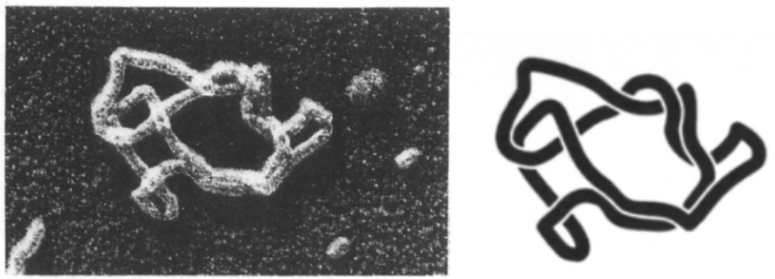}
\caption{Electron micrograph of a five-crossing knot arising from Gin recombination.  Figures from~\citet{kanaar1990processive} (permission obtained).  
}\label{fig:real.knot}
\end{figure}

Actions of site-specific recombinases can produce a wide variety of possible knots under laboratory conditions~\citep{dean1985duplex}, such as that shown in Figure~\ref{fig:real.knot}. %
Interestingly, in bacteriophage capsids the distribution of knot types in the wild is not uniform across knot types, or even across knots of the same crossing number~\citep{arsuaga2005dna}. For instance, the achiral figure-of-eight knot $4_1$ is surprisingly scarce.  
The observed distribution of knots suggests that the knot type may influence fitness, and so may carry some information about metabolism, %
or it may give information about the mechanism that gives rise to knotting. 

As noted, topoisomerase plays an indispensable  role in cell replication.
When cells replicate, the two strands in the double helical DNA split and complementary nucleotides are synthesized along each of the arms of the replication, in a process akin to unzipping a zip.  When the DNA is circular, replication begins at an \emph{origin} at which the two strands are pulled apart and replication proceeds along the forks on each side of the origin, ending at a point called the \emph{terminus}.  The problem is that this process cannot resolve the twisting that occurs without cutting the DNA at some point.  Enter topoisomerase.  Various types of topoisomerase exist, one cutting a single strand to allow untwisting and reconnecting to occur (Type I)~\citep{dean1985duplex}, another cutting both strands of the double helix, as described in Figure~\ref{fig:dna-cutting-rejoining}A and~\ref{fig:dna-cutting-rejoining}B  (Type II)~\citep{liu1980type}. %
In other words, the DNA cutting and rejoining that is effected by topoisomerases is essential to the reproductive processes of bacteria~\citep{zechiedrich1995roles,Darcy-Applications-Banach1998,alexandrov1999mechanisms,postow1999knot,postow2001topological,hardy2004disentangling}.  Consequently, understanding topoisomerase action is a goal for the development of anti-bacterial and anti-cancer drug treatments (for example,~\citet{schmittel2010controversies}).  

Knotted structures are regularly being observed in DNA, and their origins are a topic of active research (e.g.~\citet{marenduzzo2009dna}).
Examples also arise in proteins, where a Stevedore knot $6_1$ (with six crossings) has been observed~\citep{bolinger2010stevedore}.  Such a knot can be generated with a single change of crossing from the unknot, provided the unknot is suitably arranged before the crossing is changed (the mechanisms in this case are unclear).

\section{Current approaches to modelling inversions}

Inversions in bacteria are often studied with a view to phylogeny reconstruction because focussing on the inversion process avoids the obfuscating effects of horizontal transfer~\citep{Darling2008}.
That is, one may attempt to reconstruct the evolutionary history behind a set of organisms under the assumption that the only evolutionary process is inversion, yielding, at least, an approximation of the true phylogeny.  In order to do this completely, one must be able to decide, for any set of related genomes, a nearest common ancestor.  This is called the ``reversal median problem", because of the link to computer science problems known as reversals (e.g. the pancake flipping problem~\citep{gates1979bounds}).  It is typically expressed in terms of attempting to find a common ancestor that minimizes the total number of evolutionary steps (or inversions required) to each genome~\citep{caprara2003reversal}.  Equivalently, one wants to minimize the average distance to each genome.  

Of course one may use distance based methods such as neighbour-joining to construct an inversion-based phylogeny (a good survey of these approaches can be found in the book~\citet{gascuel2005mathematics}).   In any case, one defines a metric on the set of bacterial genomes, given by setting $d(X,Y)$ to be the minimal number of inversions required to transform genome $X$ to genome $Y$.  
The genome itself can be considered to be a word in the alphabet $\{A,C,G,T\}$ by fixing a starting point on the circle, but it is more natural in this context to consider instead the genome as a sequence of \emph{genes} (ignoring intervening DNA), or even better as a sequence of preserved \emph{regions} among a given set of genomes~\citep{Darling2008}.

The first modelling of inversions in this way was as a permutation of the set $\{1,2,\dots,n\}$ defined by 
$k\mapsto j+i-k$ if $i\le k\le j$ and $k\mapsto k$ otherwise.
In other words, the sequence $i,i+1,\dots,j$ is reversed.  
The initial statement of the inversion distance problem, made in~\citet{Watterson-chrom-reversal-1982}, numbered a set of gene loci common to both genomes, and considered inversions defined as above.  %
Being on a circle, they ask for the minimal number of inversions between two genomes, without regard to either mirror images or rotations around the circle. 
That is, an arrangement of regions in the order $(a,b,c,d)$ around the circular genome is the same arrangement in three dimensions as $(b,c,d,a)$, and even $(d,c,b,a)$, because the difference is merely a rotation or reflection of the whole genome.  In group-theoretic terms, one might say that these arrangement are equivalent \emph{up to the action of the dihedral group}~\citep{egrinagy2013group}.

Subsequent work, mainly by bioinformaticians and computer scientists, treated a form of the problem in which the chromosome is linear, rather than circular.  This is based on the modelling assumption that inversion events are equally likely, irrespective of the length of DNA inverted.  This model led to several interesting algorithms, usually involving a translation into a graph theory problem~\citep{kececioglu1993exact}.  
Research into the problem then shifted in two directions: to treating \emph{signed} inversions, and to finding an actual sequence of inversions that realizes the inversion distance, sometimes called \emph{sorting by reversals}~\citep{Bafna1993genome}.  

\begin{figure}
\includegraphics[width=84mm]{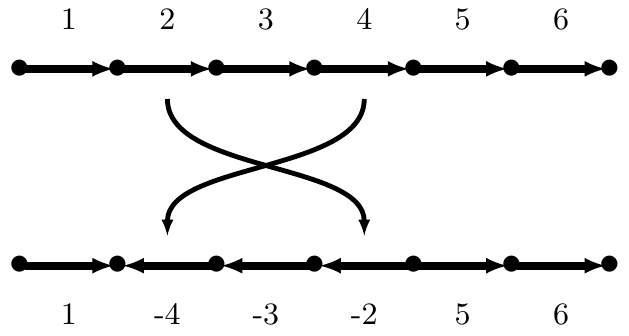}
\caption{Signed inversion on the regions $[2,3,4]$, sending them to $[-4,-3,-2]$.}\label{fig:signed}
\end{figure}

A signed inversion takes a sequence $(a_1,\dots,a_i)$ and not only reverses the order but changes the sign:  $(a_1,\dots,a_i)\mapsto(-a_i,\dots,-a_1)$, as in Figure~\ref{fig:signed}.  
Effectively, signed inversions track not only the position of the region on the genome, but also the orientation.  Giving each gene an orientation as well as a relative position is more realistic, and surprisingly turns out to be more tractable (the unsigned inversion sorting problem has been shown to be NP hard~\citep{caprara1997sorting}).  Using signed inversions is more realistic because it takes account of the fact that each preserved region has an implied orientation. There is an inherent polarity in DNA built into the biochemistry, with the two strands of the double-helix having mutually opposite orientations: inverting a segment of double-stranded DNA results in each strand of the inverted segment joining to the remaining parts of the complementary strand. %
A polynomial time algorithm for finding the minimum number of signed inversions between two genomes was presented in~\citet{hannenhalli1999transforming}, based on the ``breakpoint graph"~\citep{Bafna1993genome}, and a \emph{linear} time algorithm was given in 2001 by~\citet{bader2001linear}. 
The breakpoint graph approaches to this problem have been translated to an algebraic formalism at a similar time~\citep{meidanis2000alternative}.  
The book by~\citet{fertin2009combinatorics} provides a survey of combinatorial methods such as these, and their links to phylogeny.

Several extensions to this family of problems have been pursued~\citep{bergeron2005inversion,li2006algorithmic,hayes2007sorting}.  Broadly these proceed in two directions: considering additional mutation processes such as transpositions or block interchanges~\citep{bader2009sorting,lin2009sorting}; or adding a cost function to the length of an inversion within the distance calculation on the basis that inversion lengths are not uniformly distributed~\citep{pinter2002sorting,swidan2004sorting,sankoff2004distribution}.  
This field in general is now reaching a mature stage of development, and has become a branch of computational algorithmics, studied in many cases without reference to biological motivation.

\section{The tangle algebra approach to topological evolution}\label{s:tangles}

The use of tangle algebras to model the processes giving rise to knotting in DNA provides an excellent --- and unfortunately uncommon --- example of the application of algebra to biology.

The tangle algebra approach to knotting in DNA began with the study of Tn3 resolvase acting on unknotted DNA to produce a range of different knots in proportions that could be placed in an order that decreased exponentially.  Because the enzyme binds to the DNA at a specific site, any topological action of the enzyme on the DNA can be considered in a small three dimensional region of the cell containing the site. This motivated the use of \emph{tangles} that had been introduced by~\citet{Conway1970}.  

In order to understand how this enzyme was acting, it was assumed that the enzyme was acting in a consistent ``processive" way at the site it was bound before releasing the DNA. The distribution of knots was then inferred to reflect the different times of release of the enzyme. The model arising from this assumption had already produced testable, and verified, predictions of knot products~\citep{wasserman1985discovery}, but the tangle algebra approach made it possible to write down tangle equations that reflected the progressive repeat action of the resolvase~\citep{Ernst-TanglesDNA-MPCPS1990,ernst1996tangle,ernst1997tangle}.  %

\begin{figure}[ht]
\includegraphics[width=84mm]{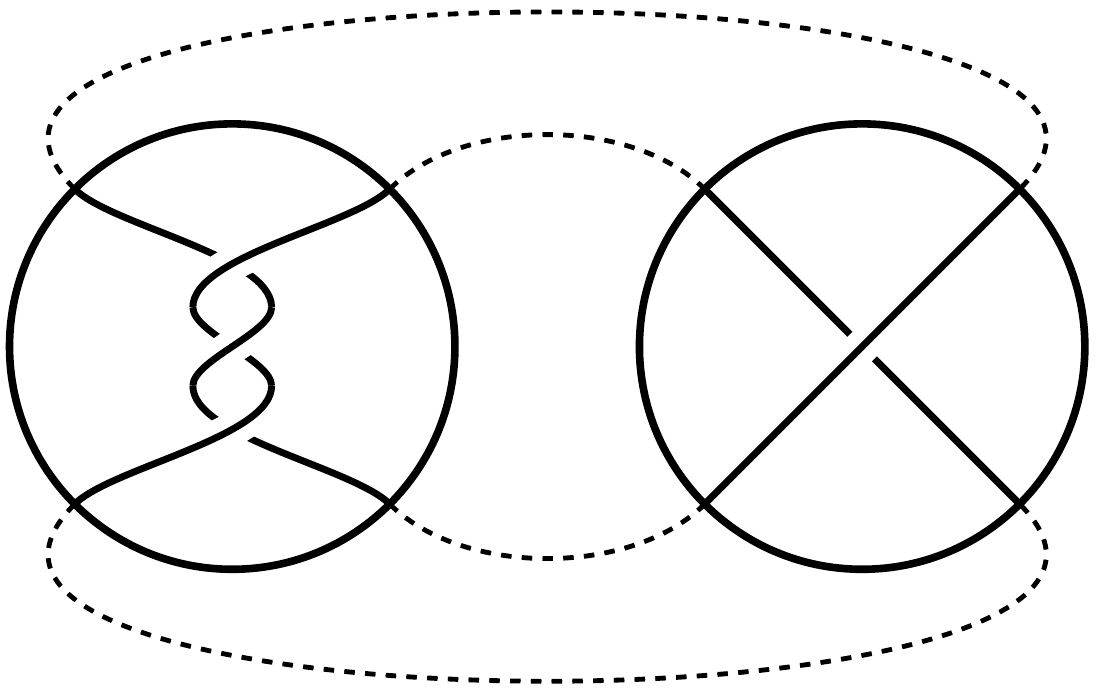}
\caption{Tangle representation of the action of Tn3.  The left tangle shows the arrangement of the substrate (before recombination the top two corners are connected, as are the bottom corners), and the addition of the right tangle gives the product of the Tn3 action, which can be repeated to generate more products. In the case of a single action as shown here, the product is a Hopf link.}\label{fig:tangle}
\end{figure}

A tangle is a box, or circle, with two strings passing through it, whose endpoints are at opposite pairs of corners of the box (consider the circle to be on the circumference of a 3-ball with strings coming from the SW, SE, NW and NE directions).  
Tangles are ``multiplied" by concatenating the boxes side by side and joining the strings up, as in Figure~\ref{fig:tangle}.  The Tn3 resolvase latches on to the synaptosome (a specific region of the DNA where strands are crossing in the right way) and through cutting and rejoining has the effect of multiplying the synaptosome tangle by another fixed tangle.  %
Under the hypothesis that the substrate was arranged in a specific way, they were able to show that these equations had a single solution, supporting the conjecture made in~\citet{wasserman1985determination} that the knotting was the result of the Tn3 staying latched over-long and acting by addition of the tangle more than once~\citep{Ernst-TanglesDNA-MPCPS1990,Ernst-Solving-ProcCambPhilSoc1999}, using the fact that tangles arising in this way are ``{rational}"~\citep{kauffman2004classification}.
Rational tangles are those that are obtained from a trivial tangle by successive twists swapping either the NE/SE strings or the SW/SE strings.
A similar approach was taken to study the effect of Gin recombinase~\citep{Vasquez-Tangle-ProcCambPhilSoc2004}.
Surveys of this and related approaches are widespread, but some good sources are~\citet{sumners1995analysis,goldman1997rational,Darcy-Applications-Banach1998,murasugi2007knot,kauffman2009tangles} and~\citet{Buck2009dna}. 

Tangles continue to be used to describe the synaptic structure giving rise to recombination, for example extending the model to include 3-string tangles~\citep{emert2000nstring,darcy2009tangle,kim2009topological,cabrera2010braid}, 
and have also been used to make predictions of the possible knots that may arise under different hypotheses about the substrate arrangement~\citep{valencia2011predicting}.  Other algebraically related approaches to DNA configurations have been attempted, including the use of the action of affine Lie groups through a translation of knot surgery into tilings of the plane on which the Lie group acts~\citep{bodner2012affine}, as well as the study of self-assembly of DNA polyhedra~\citep{hu2011topological}.

\section{An algebraic view of local evolutionary processes}
\label{sec:alg.inversions}

Within the \emph{local} view, the genome is fundamentally a circular sequence of nucleotides, represented by the alphabet $\mathcal A=\{A,C,G,T\}$, and so can be thought of as an element of a formal language  $\A^n$ for some integer $n$ in a realistic range for the size of a genome of between $10^5$ and about $10^7$.  Of course, this should be considered modulo the action of the group of symmetries that acts on the genome, thought of as a circle with $n$ evenly-spaced points representing the nucleotides.  In the case of a genome configured as an unknot (a simple circle), the group is the dihedral group on $n$ letters, $\D_{n}$.  

In the context of studying relationships among a set of genomes, and as noted in the Introduction, it is not always helpful to make comparisons at the nucleotide level because the resolution is too fine.  Instead, one can look at the set of genes on the genomes, and more recently, use genome sequence data to identify regions of DNA (effectively words in the alphabet $\A$) that are preserved (up to orientation) in each of the genomes under investigation~\citep{darling2004mauve}.  This may seem counter-intuitive, since the goal would seem to be to identify points of difference among genomes rather than similarity.  However, for investigating the action of mutations such as inversion, one seeks to ignore base-pair changes (at nucleotide level) to focus on the movements of larger segments.   Recent studies of sets of eight genomes of \emph{Yersinia pestis}, the cause of bubonic plague, have found between 60 and 80 preserved regions~\citep{Darling2008,liang2010genome}.  In one of these studies, for example, each genome can be considered to be a permutation of the 78 regions $\{R_1,\dots,R_{78}\}$ modulo again the dihedral group, in this case $\D_{78}$. 

Given a set of genomes related by inversions, and
a set of $n$ regions of DNA that are common to each genome, inversions can be thought of as generators of a group acting on the set of possible genomes that permutes these $n$ regions.  The set of all possible inversions relating these genomes is then the set of signed permutations on these regions, which is isomorphic to the finite Coxeter group $W$ of type $B$, also called the hyperoctahedral group.  One may think of the elements of this group as signed permutation matrices, in which each row and each column contains exactly one non-zero entry, which is either $1$ or $-1$.  The finite Coxeter groups are groups generated by reflections in real Euclidean space, and are well-studied, having connections to many parts of algebra.  In particular, they arise in the representation theory of finite groups of Lie type, where they appear as  Weyl groups (see for example~\citet{Hum90}).  In that context they are studied with respect to a particular presentation (generating set and relations) that corresponds with the Dynkin diagrams that arise in Lie theory.  The standard generators for the type $B$ Coxeter group are the transpositions $s_{i}=(i\ i+1)$ and the map $t$ sending $1\leftrightarrow -1$ and fixing all $i>1$.
In this framework, we regard the set of inversions as a subgroup of $W$ generated by the biologically plausible signed permutations that the model allows.

Since such inversions generate the whole group, any pair of genomes are connected by a unique group element, which may potentially be represented by a number of different sequences of inversions.  The minimal number of inversions required to write this group element we may call the \emph{inversion length} of the group element, and is precisely the inversion distance between the genomes.  
In other words, the inversion distance problem is translated into the question of the behaviour of a {length function} with respect to a set of non-standard generators representing the inversions.  It has been known for a long time that in general, given a set of generators for a permutation group, finding the length of a given permutation in terms of those generators is NP hard~\citep{even1981minimum}. 
Clearly this is not the case for all \emph{particular} groups or sets of generators: for instance, using standard generators for a finite Coxeter group~\citep{Hum90}.  At issue in the context of inversions is whether it is possible to find the length of a group element in terms of the generators mandated by the biological model. 

The algebraic approach provides an alternative framework that can be generalized to more sophisticated and realistic scenarios.
For example, it is known that inversions fix the position of the {terminus} of replication relative to that of the origin, breaking the genome into two evenly balanced replichores~\citep{eisen2000evidence}. %
While this fact has been incorporated in a limited way into recent approaches to the inversion distance problem, a group-theoretic framework makes this restriction simple to represent as the stabilizing subgroup of the terminus (existing approaches involving fixing the terminus have assumed the inversions are symmetric about the origin~\citep{ohlebusch2005median}).  Group theorists have studied alternative length functions on these groups, and it is possible progress will be made along similar lines~\citep{dyer1990reflection,howlett1999reflection}.  
Taking a group-theoretic approach allows the translation of related questions, such as that of reconstructing phylogeny, into algebraic questions, for instance about the Cayley graph of a group~\citep{moulton2011butterfly}.  

Recent work by the author and collaborators has applied this group theoretic approach to the inversion distance problem for a model in which inversions act on only two regions at a time, and in which orientation is ignored~\citep{egrinagy2013group}.  When the number of regions that inversions may act on is restricted, many standard approaches fail because they assume any inversion of a region is as likely as the inversion of the complementary region, allowing the problem to be considered as if the genome were linear.  That is, any model in which inversions act on only a restricted number of regions must account for the circular structure of the genome.  This kind of model is not unrealistic because it is also known that inversions of shorter segments occur more frequently than those of long segments~\citep{Darling2008}.    In~\citet{egrinagy2013group}, the need to work with circular permutations is handled by lifting the permutation to the group of periodic permutations of the integers, known as the \emph{affine symmetric group}.  Results about the length function in this group are then able to be used (in particular a length formula given by~\cite{shi1986kazhdan}).  While this lifting from circular to affine permutations is not trivial and provides some theoretical challenges, it nevertheless produces a polynomial time algorithm and indicates that the affine symmetric group is the ``right" place to study circular permutations.  This is especially the case when the assumption of uniformly likely inversions is removed from the model.

The application of group theory and other algebraic ideas to local evolutionary mechanisms will also allow generalizations such as the incorporation of other known types of mutation into this model. %
For instance, translocation is another invertible operation that can be studied. %
Deletion requires more care as it is not invertible: while horizontal transfer does involve insertion, and so could be considered an inverse operation, it generally does not reinsert a piece recently deleted from the same chromosome.  Incorporating non-invertible actions will require modelling the action as that of a semigroup rather than a group. %
For all of these extensions and variations, extensive theoretical and computational tools that have been developed within the world of algebraic research can be brought to bear.  For instance, the power of computational systems such as GAP~\citep{GAP4} and Magma~\citep{bosma1997magma} (and their many packages) are barely used in biology (the aforementioned work of the author implements the algorithm using GAP~\citep{egrinagy2013group}).

Finally, it is worth noting that many assumptions made about inversions, such as the fixing of the terminus of replication, are actually not quite so rigid.  It is more correct, for instance, to say that the terminus stays within a small distance of the antipode of the origin.  Some statistical approaches have already been taken to the inversion distance problem (e.g.~\citet{york2002bayesian}, \citet{miklos2005genome}), and the logical development of this theme is to use models involving group actions in a probabilistic setting --- a genuinely multidisciplinary endeavour.  Integrated computational systems such as Sage~\citep{sage}, that can call on GAP or Magma as well as statistical packages such as R~\citep{R2011}, are likely to play an important role.

\section{Inversions and knotting in a common modelling framework}\label{sec:common.framework}

\begin{figure}[ht]
\includegraphics[width=84mm]{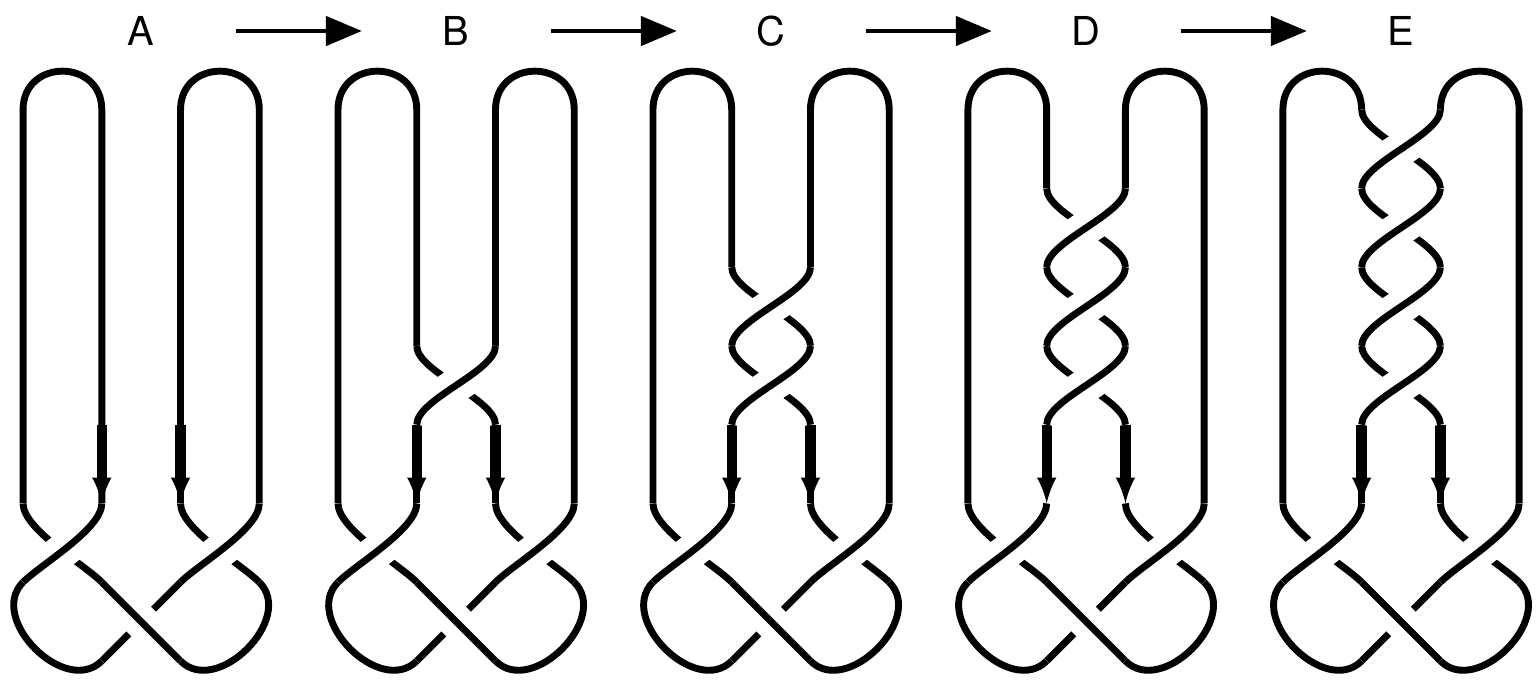}
\caption{The action of resolvase enzymes Tn3 and $\gamma\delta$ on the \emph{res} site on a substrate (A) that is twisted gives rise to a series of products through processive recombination.  The products are (B) a Hopf link (catenane), (C) a figure-8 knot in which the sequence is as in the original substrate, (D) another link (the Whitehead link), 
and (E) a six-crossing knot.  Figure adapted from~\citet{grindley2002movement}.} %
\label{fig:knot.process}
\end{figure}

\begin{figure}[ht]
\includegraphics[width=84mm]{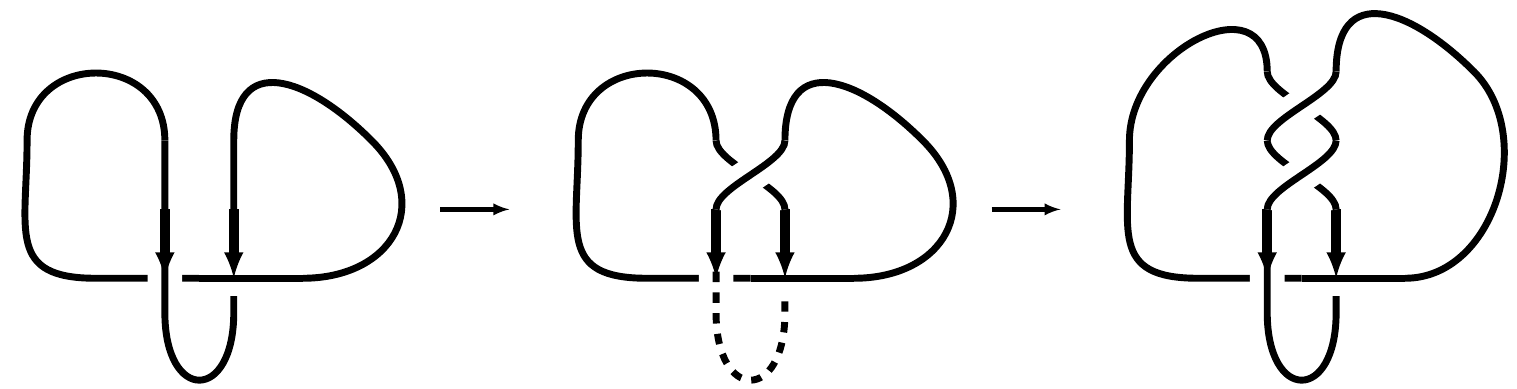}
\caption{Schema of the Gin inversion synapse, as in~\citet{klippel1993analysis}.  After the first recombination event, the segment below the  \emph{gix} recombination sites, shown dashed, is inverted relative to the rest of the DNA.  After the second event, this orientation is restored but the DNA is knotted in a trefoil configuration.}\label{fig:inversion.process.lit}
\end{figure}

\begin{figure}[ht]
\includegraphics[width=84mm]{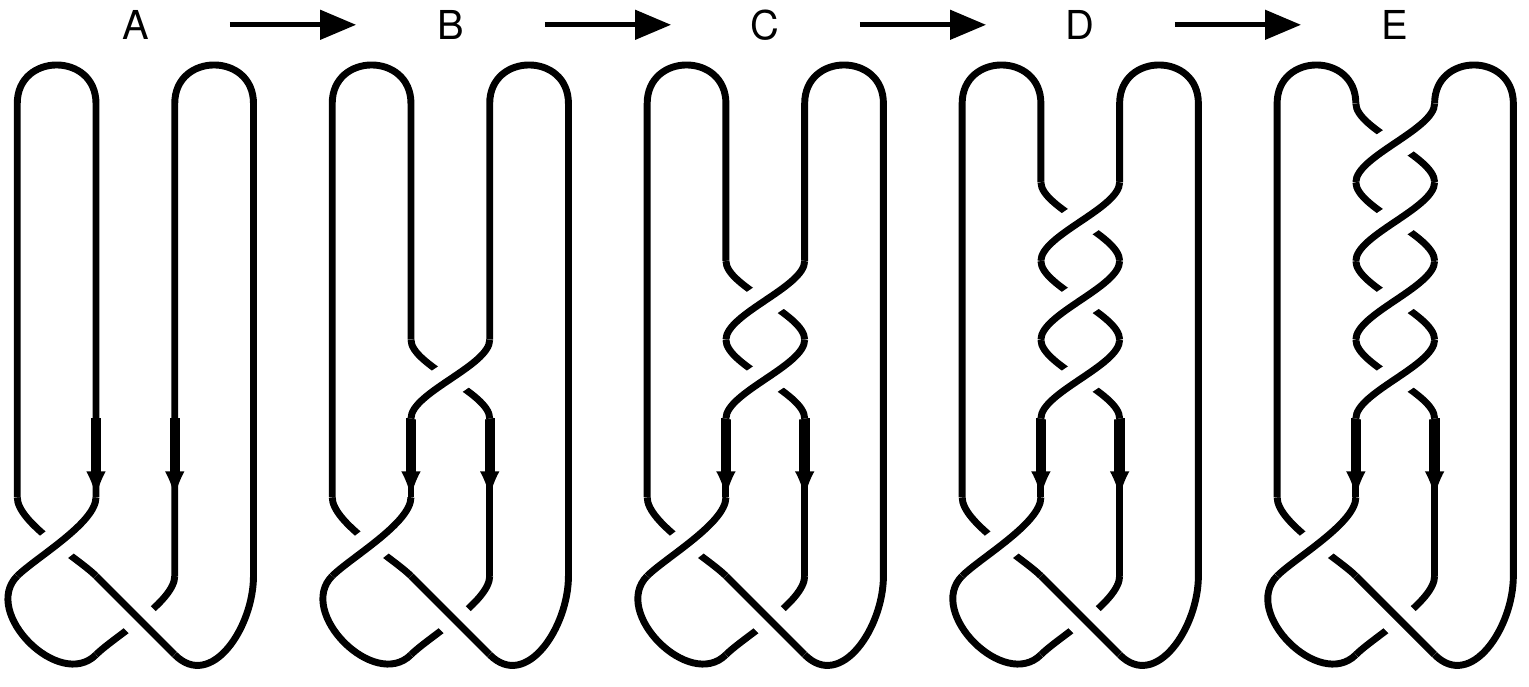}
\caption{An alternative representation of the inversion process, adapted from Figure~\ref{fig:inversion.process.lit}.  Note that rearrangements A, C and E all have identical sequence data, whereas B and D are inversions.  This illustrates the interplay between inversion and knotting, because while A, C and E are isomers, their knotting is distinct, giving the trivial knot, the trefoil, and a six-crossing knot respectively.  The inversions B and D themselves are not topologically identical, with B being the trivial knot and D a five-crossing knot.} %
\label{fig:inversion.process}
\end{figure}

While the questions about inversions and about knotting that we have described above tend to be addressed separately in the modelling literature, the biological mechanism giving rise to inversions and knotting is widely acknowledged as being the same: cutting and rejoining of DNA double-helices.  In the case of the action of site-specific recombinases 
the difference between knotting and inversion is in the arrangement of the DNA (the substrate) when the recombination takes place.   For instance, the action of the resolvase enzymes Tn3 or $\gamma\delta$ transposons at the \emph{res} site give rise to knotting and linking when the sites are aligned in the same orientation (as direct repeats) and the DNA is twisted in a certain way~\citep{li2005structure,grindley2006mechanisms,jayaram2009difference} (see Figure~\ref{fig:knot.process}).  (Recall from Section~\ref{s:tangles} that Tn3 is an enzyme whose action has been analyzed using tangle algebras). 
On the other hand, the action of the Gin recombinase gives rise to inversion when the two \emph{gix} sites  
are aligned with the same orientation and the DNA twisted slightly differently (see Figures~\ref{fig:inversion.process.lit} and~\ref{fig:inversion.process})~\citep{klippel1993analysis,grindley2002movement,jayaram2009difference}.

\begin{figure}[ht]
\begin{center}
\includegraphics[width=84mm]{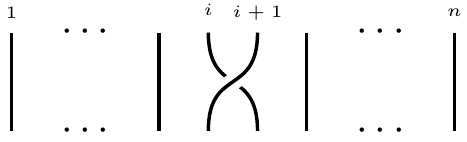}
\end{center}
\caption{The braid group generator $\sigma_i$.}\label{fig:sigma_i}
\end{figure}

The recombination events represented in Figures~\ref{fig:knot.process} and~\ref{fig:inversion.process} can be represented as braid closures in a remarkably elegant way.   In each case there is a substrate braid (``plat-closed''), and the action of the recombinase is to add a generator $\sigma_i$ (see Figure~\ref{fig:sigma_i}) to the substrate braid before it is plat-closed~\citep{sumners1995analysis}.

\begin{figure}[ht]
\includegraphics[width=84mm]{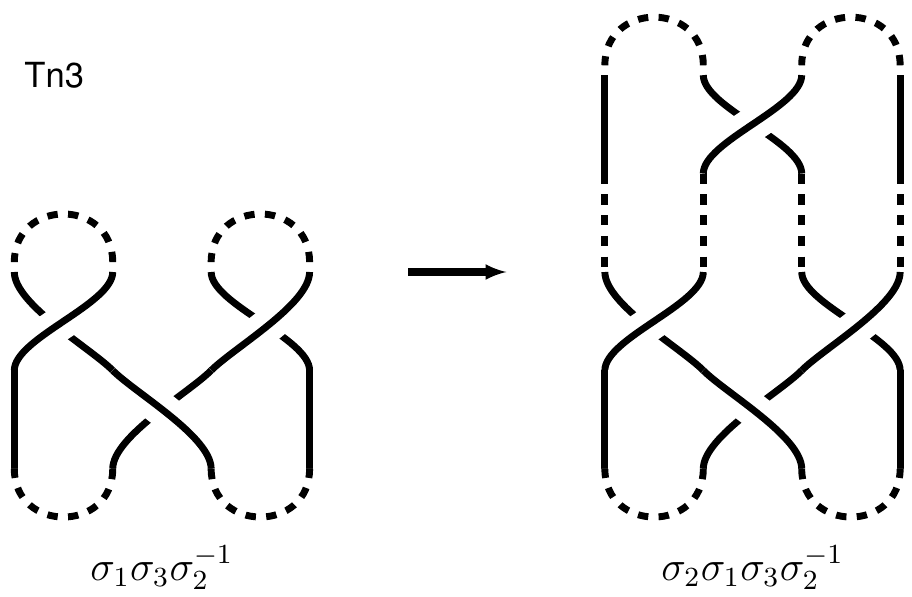}
\caption{The action of Tn3 shown as the plat closure of a braid (see Figure~\ref{fig:inversion.process.lit}, and Figure~\ref{fig:tangle} for the tangle version).  }\label{fig:Tn3.plat}
\end{figure}

A braid on $n$ strings is a set of $n$ strings joining two parallel lines of $n$ points, such that the strings pass continuously downwards.  The set of braids on $n$ strings forms a group whose multiplication is performed by placing one diagram below the other and joining corresponding strings (an example is in Figure~\ref{fig:Tn3.plat}). The \emph{braid group} is generated by braids in which an adjacent pair of strings is interchanged.  We denote these generators by $\sigma_i$, being the braid that interchanges strings $i$ and $i+1$, with the $i$'th string passing behind (Figure~\ref{fig:sigma_i}).  The inverse $\sigma_i^{-1}$ is the same but with the crossing reversed so that the $i$'th string comes in front; in general the inverse of any braid can be obtained by taking its reflection in a mirror placed below the lower $n$ points.  This is easy to see in the context of the operation in the group being stacking diagrams on top of each other.  Good references for an introduction to braid groups are the books by~\citet{murasugi1999study} and~\citet{kassel2008braid}.

Braids are widely used in the mathematical study of knots.  There are two standard ways in which a braid can be transformed into a knot (or link).  The first, more common way, is to join the string in $i$'th position at the bottom with the string in the same position at the top. The second way is possible with an even number of strings, and involves joining adjacent pairs of strings at the top and the bottom.  The latter is called \emph{plat closure}.  An example of the plat closure of a braid is given by the dashed lines at top and bottom of Figure~\ref{fig:Tn3.plat}.

In the case of the action of Tn3 (shown in Figure~\ref{fig:knot.process}), giving rise to knotting, the base braid is $\sigma_1\sigma_3\sigma_2^{-1}$, and the recombinase acts by adding a power of $\sigma_2$ as a prefix.
The first step of this process is shown in Figure~\ref{fig:Tn3.plat}.  The family of products of processive recombination of the action are plat closures of the braids 
\[\sigma_2^i(\sigma_1\sigma_3\sigma_2^{-1}),\] 
where $i$ is the number of twists added due to the recombinase staying bound to the substrate.  
On the other hand, in the case of the Gin inversion synapse (Figure~\ref{fig:inversion.process}), we see the substrate has one fewer twist and the base braid is $\sigma_1\sigma_2^{-1}$.  The sequence of products of processive recombination is then given by the plat closure of
\[\sigma_2^i(\sigma_1\sigma_2^{-1})\] 
for $i\ge 1$.

Note, given a certain substrate, the form of the braid whose closure gives the substrate is not unique.  For instance, in the substrate in Figure~\ref{fig:knot.process} we could push all crossings  except the last to the left pair of strands, so that the substrate could instead be produced from the plat closure of $\sigma_1^2\sigma_2^{-1}$, instead of $\sigma_1\sigma_3\sigma_2^{-1}$.  The equivalence of braids under the standard closure is well-studied (braids produce the same knot under standard closure if they can be reached from each other via a sequence of Markov moves~\citep{markov1935freie,birman1974braids}), and similar information is known for plat closure~\citep{birman1976stable,tawn2008presentation}.  This raises the question of whether more actions of recombinase are performed that are not detected by changes in topology or sequence, and hence go undetected by experiment.

Despite this, we do know that every knot can be represented as the plat closure of a braid, and hence that \emph{every} substrate for a recombination reaction can indeed be represented in this way.  (To see that any knot can be represented as the plat closure of a braid, imagine placing the knot between two parallel lines, and pulling loops that turn up down to the bottom line, and loops that turn down up to the top line). 
In addition, we are able to derive some standard results on the path of the processive recombination when the substrate 
can be arranged as in Figure~\ref{fig:inversion.process}(A). %
For a start, it is always possible to draw such a substrate as the plat closure of a braid on \emph{four} strands, known as a \emph{4-plat}.  Second, the braid word $\s_1^{k-1}\s_2^{-1}$ can be used to express the simplest case of an unknotted substrate with $k$ crossings.  The study of recombination events on substrates that can be arranged as 4-plats was part of the motivation for the tangle model of site-specific recombination due to~\citet{Ernst-TanglesDNA-MPCPS1990}.

Consider an arbitrary orientation on the unknotted chromosome, arranged as in Figure~\ref{fig:inversion.process}(A).  If there are an even number of crossings, then when arranged so that all crossings are below the recombination sites, the middle two strands (where recombination occurs) will have opposite orientations (the recombination sites are in inverted repeat). 
Then recombination has the effect of joining strands of opposite orientation, resulting in an inversion.  Subsequent recombinations alternate between re-orienting and inverting.  In terms of braids, if we represent the substrate with $2k$ crossings as $\sigma_1^{2k-1}\sigma_2^{-1}$, we have that the products $\sigma_2^i(\sigma_1^{2k-1}\sigma_2^{-1})$ are inversions when $i$ is odd and preserve the orientation when $i$ is even. 
On the other hand if the substrate is arranged with an odd number of crossings, we may write the braid in form  $\sigma_1^{2k}\sigma_2^{-1}$ for some $k$.  In this case if we follow the strand in the second position from the left down through the braid and up the other side after plat closure, it emerges on the outside (fourth) strand.  The effect of the additional twist given by recombination is to rejoin this strand with the second strand, and the result is a link.  Subsequent events alternately restore the strand to a single loop or create links, so that the product $\sigma_2^i(\sigma_1^{2k}\sigma_2^{-1})$ is a link precisely when $i$ is odd.

\begin{figure}[ht]
\begin{center}
\includegraphics[width=175mm]{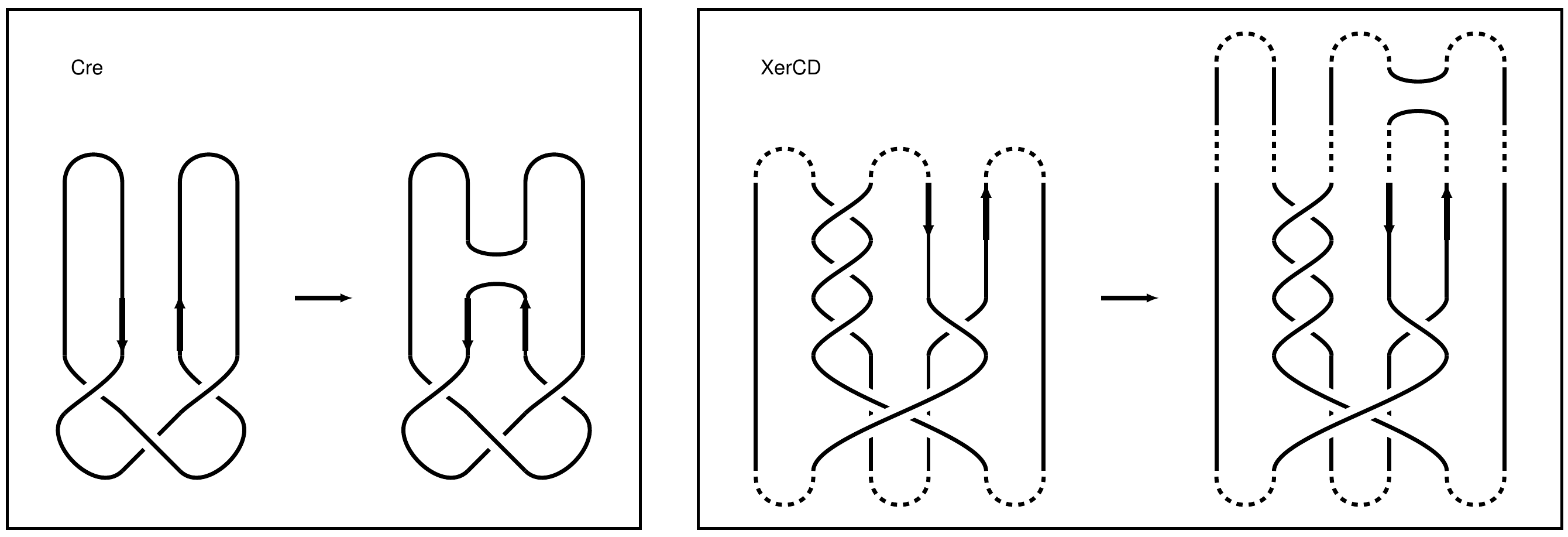} %
\caption{On the left, Cre recombination when the sites are inverted repeats, giving rise to a trefoil knot (adapted from~\citet[Figure 2]{kim2009topological}; see also \citet{vetcher2006dna}). 
On the right, a model of synapsis for  XerCD-mediated recombination proposed in~\citet{reijns2005mutagenesis}. The pictorial representation of the process is adapted to illustrate how the change given in Figure~\ref{fig:transf}(C) 
can be represented as the plat closure of a BMW tangle diagram. 
}\label{fig:XerCD}
\end{center}
\end{figure}

\begin{figure}
\includegraphics[width=84mm]{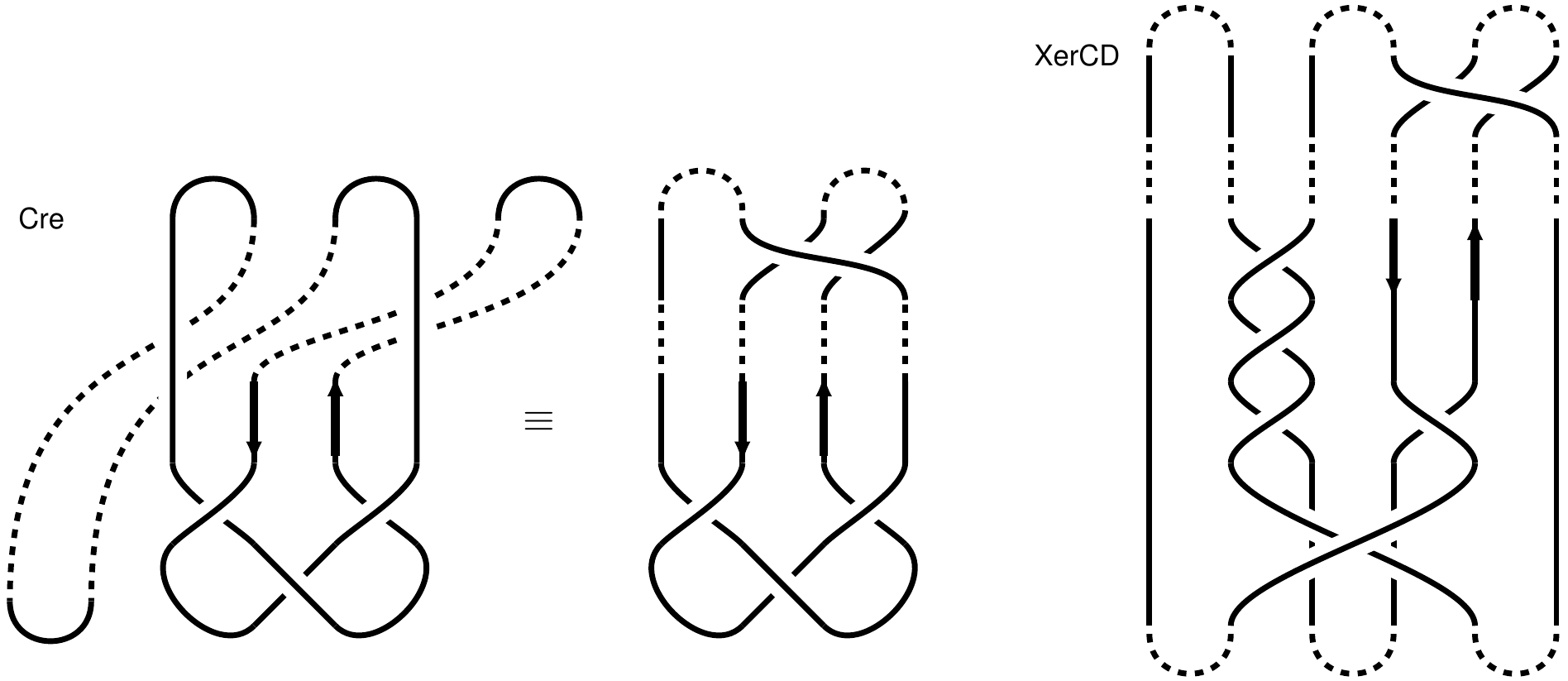}
\caption{Cre and XerCD recombination shown as the plat closures of braids.  Contrast with the representation in Figure~\ref{fig:XerCD}, in which the same process is represented using BMW algebra generators.  Note, there are many alternative ways one might represent these processes as plat closures.
}\label{fig:XerCD.braids}
\end{figure}

The recombinations described in the previous paragraph involve transformations of form given in Figure~\ref{fig:transf}(A) or ~\ref{fig:transf}(B).  However, not all recombination events produce rejoining in this form.  Some, such as the action of Cre on the site \emph{loxP}~\citep{vetcher2006dna,kim2009topological}, or the action of XerCD to resolve plasmid dimers~\citep{reijns2005mutagenesis}, have the effect shown in Figure~\ref{fig:transf}(C).  These recombinations are shown here in Figure~\ref{fig:XerCD}. 
While these reactions initially appear not amenable to the braid analysis because the strings turn upwards, it is possible for the substrate, and even the action, to be represented as a braid, albeit with a bit more effort.  This is because even if the diagram is drawn with upturned strands, these up-turnings %
(respectively down) can be pulled down (resp. up) into the plat closure, as shown in Figure~\ref{fig:XerCD.braids}.  (It should also be noted that the three-dimensional synaptic arrangement for the  same reaction can often be projected onto the plane with apparently different alignments of recombination sites~\cite{sumners1995analysis,vazquez2005tangle}).

\begin{figure}[ht]
\begin{center}
\includegraphics[width=84mm]{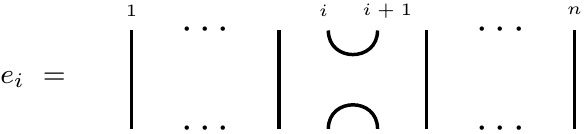}
\end{center}
\caption{The BMW algebra generator $e_i$.}\label{fig:e_i}
\end{figure}

While braid groups can be used to model such actions, there are other diagram algebras that may provide an alternative model (an \emph{algebra} is a set with two operations consisting of linear combinations of elements that can be multiplied --- like a vector space in which we can multiply).  The Birman-Murakami-Wenzl algebra (or BMW algebra)~\citep{birman1989braids,murakami1987kauffman} is similar to braids in that its basis elements consist of strands that connect two lines of $n$ parallel points.  However in this algebra the strings are allowed to return up (or down) to the level of their origin. %
The set of such diagrams no longer forms a group because not all diagrams have ``inverse" diagrams.  However a multiplication can be defined in the same way as braids, and linear combinations of such diagrams form an algebra over Laurent polynomials in two variables $\Z[\lambda_1^{\pm 1},\lambda_2^{\pm 1}]$ (see~\citet{birman1989braids,murakami1987kauffman}  for more details).  The generators of the BMW algebra include the braid group generators $\sigma_i$ as shown in Figure~\ref{fig:sigma_i} as well as the generators $e_i$ given by the diagrams shown in Figure~\ref{fig:e_i}.  

\begin{figure}[ht]
\begin{center}
\includegraphics[width=150mm]{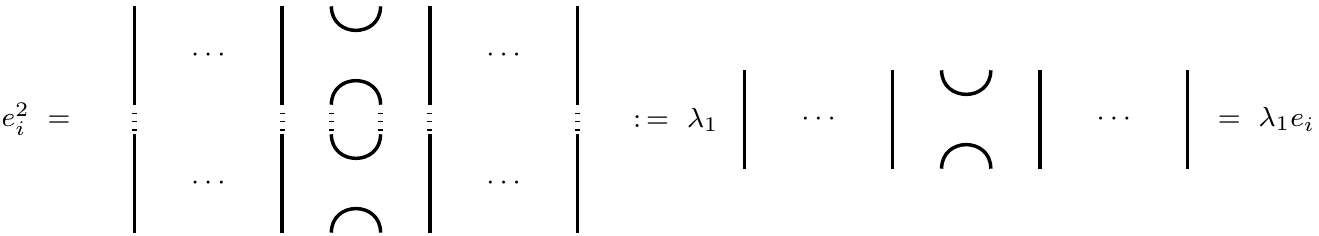}
\end{center}
\caption{The BMW algebra relation $e_i^2=\lambda_1 e_i$.}\label{fig:bmwrel}
\end{figure}

As a consequence of its origins arising from knot invariants, the multiplicative structure of the BMW algebra is not given simply by concatenation of diagrams, but includes some additional relations such as $e_i^2=\lambda_1 e_i$ and $\s_ie_i=e_i\s_i=\lambda_2 e_i$, for some parameters $\lambda_i$, to account for closed loops and twists arising in products (see Figure~\ref{fig:bmwrel} for an illustration of the first of these). 
There is also a ``skein relation" $\s_i+\s_i^{-1}=\lambda_3(1+e_i)$ that arises out of the requirement for the algebra to have a trace that relates to Jones' knot invariant~\citep{Jones1985polynomial}. %

The substrates of Cre and XerCD shown in Figure~\ref{fig:XerCD} can be expressed in braid notation as plat closures.  For instance the substrate of XerCD can be written as a product of the braid group generators as 
$\s_2^3\s_2^{-1}\s_4^{-1}\s_4\s_3\s_2\s_4^{-1}$ (which may be simplified to $\s_2^{2}\s_{3}\s_{2}\s_{4}^{-1}$).  The recombination action may then be written as an element of the BMW algebra as the plat closure of $e_{4}(\s_2^3\s_2^{-1}\s_4^{-1}\s_4\s_3\s_2\s_4^{-1})$.
This product could also be represented in terms of braids alone by the plat closure of the six-strand braid $\s_{4}^{-1}\s_{5}^{-1}(\s_2^3\s_2^{-1}\s_4^{-1}\s_4\s_3\s_2\s_4^{-1})$, as shown in the right hand picture in Figure~\ref{fig:XerCD.braids}.  

In both Cre and XerCD the representation of the recombination as a plat closure of braids instead of in the BMW algebra incurs a minor penalty, namely replacing multiplication by a single BMW algebra generator ($e_2$ or $e_4$) with multiplication by a product of two braid group generators.  This is a general property: multiplication by $e_i$ (for $i$ even) and multiplication by $s_is_{i+1}$ (or $s_i^{-1}s_{i+1}^{-1}$) are equivalent as plat closures.  But it is also somewhat special because it depends on the action of the recombinase being expressible with the action of $e_i$ at the top.  While one of the features of the tangle algebra approach is that it is able to cover a wider variety of circumstances, it is also worth  considering whether the connections provided by other algebraic structures may be sufficient compensation for a possibly restricted sphere of operation.
The moral of the story is that there may be a range of possible formalisms using diagram algebras that preserve the synaptic structure observed through experiment.  In addition to the insights afforded by the tangle algebra approach, the above discussion demonstrates that these recombination processes may also be viewed as transformation in the braid group or in the Birman-Murakami-Wenzl algebra.  These alternative algebraic approaches may lead to connections with algebraic models of local evolutionary processes such as inversion, as described in the next section.

\section{Algebraic connections}\label{s:alg.conn}

We have described a modelling framework using plat closures of braids that can incorporate the biological mechanisms behind numerous recombinase actions giving rise to inversion and knotting.  However there is an important piece of the inversion story that is omitted by this analysis, namely the role of selection, at least in the case of bacterial evolution.  Selection is presumably behind a number of features of bacterial genomes that have been observed, and that we have described above (Section~\ref{sec:alg.inversions}). 
For instance, the observation that the terminus of replication is always close to the antipode of the origin is a consequence of a fitness cost due to the effect on the replication process of unbalanced replichores~\citep{eisen2000evidence}.  The models of inversion using braids or tangles do not (and probably cannot) take into account the relative positions of the origin and terminus.  Similarly, there is also evidence, described earlier, that shorter inversions occur more frequently than longer inversions~\citep{Darling2008}, and this is not a feature that is immediately accessible through a braid or tangle analysis, because braids and tangles are preserved under isotopy while specific locations on strands are all equivalent under isotopy.  Both of these selection features can, however, be studied using finite reflection groups such as the Coxeter group of type $B$, as described in Section~\ref{sec:alg.inversions}. %

Limitations are inevitable in any model: models cannot address every question that may arise.  In the evolutionary problems we have described there are two approaches that have links with algebra that either have already, or have the potential to provide biological insight: the tangle and plat closure models of the topological effects of DNA recombination; and the Coxeter group models of inversions.  The interesting point is that these algebraic structures have links that may yet bring forth a family of models that unifies both approaches.

We will now describe some of the algebraic connections that exist between the hyperoctahedral group, or Coxeter group of type $B$, and the tangle algebras.
These algebraic connections are at present not used in any biological context; the purpose of outlining these connections is in the hope that they might provide the key to a unified picture.
The first algebraic connection to detail is that between the Coxeter group of type $B$ and the braid groups.

\subsection{Coxeter groups and braid groups}
\begin{figure}[ht]
\includegraphics[width=160mm]{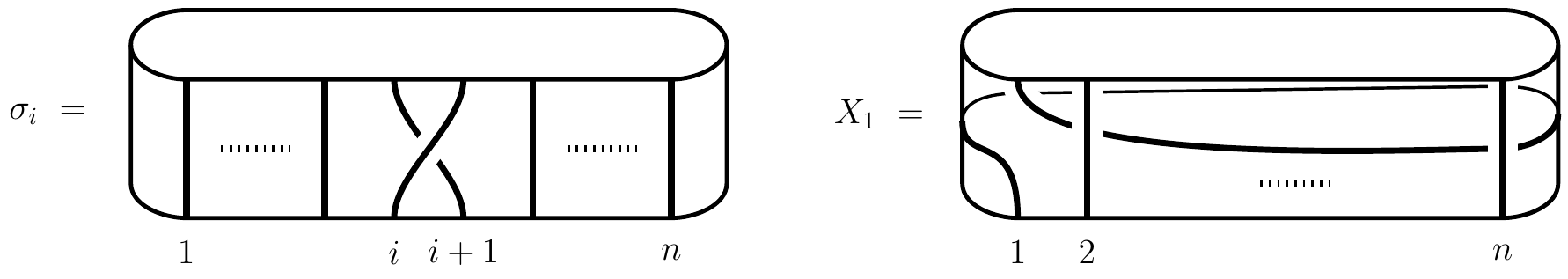}
\caption{The generators of the affine braid group.}\label{fig:affine.braids}
\end{figure}

The Coxeter group of type $B$ can be looked at in many useful ways, one of which is as a quotient of the \emph{affine} braid group $\widetilde{Br}_{n}$.  As we have seen, a braid on $n$ strings is a set of $n$ intertwined strings falling from one set of $n$ points on a line to another set on a lower line.  Affine braids have the sets of points  at the top and at the bottom of the rim of a hollow \emph{cylinder} (with slightly thickened walls), or equivalently can be thought of as \emph{periodic}. (An alternative view has them as regular braids with a rigid pole at one end that strings can loop around -- see for example~\cite{orellana2004affine}).  The affine braid group is generated by the usual generators of the braid group $\s_{i}$ for $i=1,\dots,n$, ($\s_n$ swaps $n$ and 1 around the back of the cylinder), as well as the braid $X_{1}$ that begins at position 1, passes behind the other $n-1$ vertical strings and returns to position 1 at the bottom.  These generators are shown in Figure~\ref{fig:affine.braids}.  

The regular braid group $Br_{n}$ can be obtained from the affine braid group $\widetilde{Br}_{n}$ by a projection that amounts to squashing the cylinder flat from the front.  In this projection $X_{1}\mapsto 1$, and $\s_{n}$ becomes redundant as a generator, swapping strands $1$ and $n$ behind the other strands.  
By lifting our view to the affine braid group we can see both the symmetric group (the Coxeter group of type $A_n$), and the Coxeter group of type $B_n$ as quotients of the same object.  A common view of the symmetric group is as a quotient of the regular braid group $Br_{n}$ by the relations $\s_{i}^{2}=1$, but beginning from the affine braid group the quotient is by these relations together with $X_{1}=1$. Taking the quotient from $\widetilde{Br}_{n}$ instead by the relations $\s_{i}^{2}=1$ and $X_{1}^{2}=1$ gives the type $B$ Coxeter group that arises in the study of inversions: the image of $X_{1}$ acts as the permutation $1\leftrightarrow -1$.

This standard approach of obtaining the Coxeter group of type $B$ from the affine braid group makes plain the contrasting viewpoints of local and topological evolution.  Coming from the braid group, the Coxeter group of type $B$ is a (signed) permutation group on the $n$ endpoints of the strands of the braids.  In contrast, through studying inversions the Coxeter group of type $B$ is a permutation group of linear regions along the whole chromosome, and there is no necessary correspondence between these regions and the endpoints of a braid visualization of the recombination process.  Hence we have two genuinely distinct ways in which this Coxeter group may arise through models of site-specific recombination. %

\subsection{Iwahori-Hecke algebras}

Other important algebras also arise as quotients of the braid groups, including the Iwahori-Hecke algebras.  These are deformations of group algebras of finite Coxeter groups such as the hyperoctahedral group and the symmetric group, and are obtained from the affine braid group algebra by taking the quotient by the relations $(\s_{i}-1)(\s_{i}+q)=0$ and $(X_{1}-1)(X_{1}+Q)=0$, where $q,Q$ are indeterminants.  This forces the relation $\s_{i}^{2}=q+(q-1)\s_{i}$ in the Iwahori-Hecke algebra.   It should be noted that these quotients of the affine braid group are linked to a specific presentation with generators $\s_{i}$, whereas the generators for the inversion group that are required are not the images of this quotient. 

The fact that Coxeter groups and Iwahori-Hecke algebras arise in closely related contexts is no surprise.  Most Coxeter groups arise as Weyl groups through the representation theory of finite groups of Lie type (these are the \emph{crystallographic} finite reflection groups).  They appear as double coset representatives of Borel subgroups in a group of Lie type.  In this context, Iwahori-Hecke algebras arise by taking these double cosets as a basis for an algebra, and are important in studying the way induced representations of the Lie group decompose into irreducible representations.  The parameters $q$ and $Q$ have a group-theoretic meaning in this original motivating way of looking at these algebras. 

\subsection{The Kauffman tangle algebra}
\begin{figure}[ht]
\includegraphics[width=40mm]{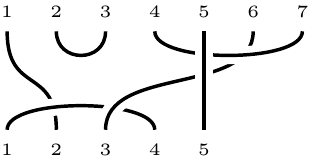}
\caption{An example of a (7,5) Kauffman tangle.}\label{fig:mn-tangle-eg}
\end{figure}

Tangle algebras themselves, aside from being seen through the prism of plat closures of braids, arise out a wider family of diagram algebras.  The Kauffman tangle algebra has a basis consisting of $(m,n)$ tangles -- diagrams of strings connecting a row of $m$ points with a row of $n$ points, in which strings may return to the level from which they originated~\citep{Morton1990knots}.  An example is shown in Figure~\ref{fig:mn-tangle-eg}. 
Tangles as used in DNA topology are generally (2,2)-tangles (there are recent exceptions, such as~\citet{kim2009topological} and~\cite{darcy2009tangle}. %
The plat closures of Kauffman tangles may be a natural way to represent recombination, because of reactions such as that of XerCD that involve a change of direction of the strand (see Figure~\ref{fig:XerCD}), and because the algebra allows flexibility in the number of nodes at each end.  In this case the action of the recombinase can be represented by multiplication by the non-braid generator $e_{i}$.  Even here, it is possible to manoeuvre the tangle to make it into a braid plat closure, as shown in Figure~\ref{fig:XerCD.braids}.  

\subsection{Knot invariants}
There is another algebraic connection through knot invariants. Tangle algebras played a key role in the development of knot invariants such as the Jones polynomial~\citep{Jones1985polynomial} and its HOMFLY generalization~\citep{freyd1985new}.  Interestingly, the use of tangle algebras in the HOMFLY paper was within a few years of their appearance in the work of Ernst and Sumners~\citep{Ernst-TanglesDNA-MPCPS1990} (they were first defined by Conway much earlier~\citep{Conway1970}).  This role of tangle algebras in knot invariants provides another link to the Iwahori-Hecke algebras, since the Hecke algebras also appear in Jones' work on his invariant.  

\subsection{Lie theory}
There is one final algebraic thread that draws these stories together.  The standard braid group $Br_{n}$ and its quotient the symmetric group are part of a Lie theory story, in that the symmetric group is a type $A$ Weyl group, arising from the representation theory of the Lie group $GL_{n}$, the general linear group.  The Iwahori-Hecke algebra arises in the decomposition of induced cuspidal representations in this case.  The group arising from inversions however is also part of the Lie theory story, because it is the Weyl group of type $B$, corresponding to the representation theory of the orthogonal group.  The algebra that plays an analogous role to that of the Iwahori-Hecke algebra for the orthogonal group is the Birman-Murakami-Wenzl algebra we have seen above.  This algebra can also be viewed as a diagram algebra of strands running between two sets of $n$ points, whose generators are $\{\s_{i},e_{i}\mid 1\le i\le n-1\}$, as discussed above in Section~\ref{sec:common.framework}.  As diagram algebras, the BMW algebra and Kauffman algebra are isomorphic~\citep{Yu-PhD-Cyclotomic-BMW-algs-2007,Morton2010basis}. %

\subsection{Diagram algebras}

Diagram algebras such as these we see here appear in widely different contexts.  While the symmetric group algebra can be thought of as the standard braids with crossings ignored, if one ``ignores crossings'' in the BMW algebra one obtains the diagrams that form the basis of the Brauer centralizer algebra, which was first defined in the 1930s and was motivated by invariant theory~\citep{brauer1937algebras}.  If one considers either the BMW algebra or the Brauer algebra and requires that strands may not cross, so that the algebra is generated by the $e_{i}$ without the $\s_{i}$, one obtains the Temperley-Lieb algebra which plays an important role in statistical mechanics~\citep{temperley1971relations}.

There is a further generalization of these diagram algebras arising from braided monoidal categories that potentially may also be applied to evolutionary mechanisms in DNA, and that is \emph{ribbon categories}~\citep{turaev2010quantum}.  These categories can be thought of as algebras whose elements consist of ribbon diagrams --- like BMW algebra diagrams in which the strings are slightly thickened --- and on which one may have ``coupons", represented by boxes on the ribbons.  In the category theory context, the ribbon bands are objects and the coupons morphisms.  In the evolutionary context the use of ribbons may enable the representation of twists on strands, while coupons might be used to represent recombination events.  This potential model remains to be explored.

\section{Concluding remarks}

The algebraic models of tangles and inversions that we have described above appear to be distinct %
algebraic stories.  These stories are part of the same big picture sitting inside various instantiations of Lie theory, and yet the way the algebraic structures arise is very different.  The algebraic models of knotting focus on the topology of strands of DNA that are equivalent up to isotopy, whereas the algebraic models of inversion focus on patterns along the strands that ignore isotopy.  And while the algebraic stories have many differences, there are several clear connections.  For instance, both the Coxeter group of type $B$ and the tangle algebras are closely connected to the affine braid group.  And while tangles arise in a type $A$ context, many site-specific recombination events could be looked at as plat closures of elements of the Birman-Murakami-Wenzl algebra --- an algebra that plays an important role in type $B$ representation theory.

Another way of viewing the key challenge in unifying these pictures is that one is a way to study inversions using braids, the other with (signed) permutation groups.  It so happens that there is an intimate link between braids and permutation groups, described above, and this gives some hope that a unified picture may be possible.
The presence of structures such as tangle algebras and Coxeter group actions in the evolutionary processes of bacterial DNA strongly suggests that there should also be a role for these related algebras.  Determining this role --- represented at the bottom of Figure~\ref{fig:phylo} --- is a key open problem requiring the expertise of algebraists working together with evolutonary biologists.

\begin{figure}[ht]
{\footnotesize  
\begin{center}
\begin{tikzpicture}[scale=0.8]
  \draw (5,6) node {phylogenetics};
  \draw (1.5,4) node {local evolution};
  \draw (1.5,2) node {Coxeter group actions};
  \draw (5,0) node {BMW algebras};
  \draw (5,-.6) node {Iwahori--Hecke algebras};
  \draw (5,-1.2) node {Affine braid groups};
  \draw[->,>=latex] (4.5,.5) -- (2.5,1.5);
  \draw[->,>=latex] (2,2.5) -- (2,3.5);
  \draw[->,>=latex] (2.5,4.5) -- (4.5,5.5);
  \draw (8.5,4) node {topological evolution};
  \draw (8.5,2) node {tangles, plat closed braids};
  \draw[->,>=latex] (5.5,.5) -- (7.5,1.5);
  \draw[->,>=latex] (8,2.5)--(8,3.5);
  \draw[->,>=latex] (7.5,4.5)--(5.5,5.5);
  \draw[<->,dashed,>=latex] (3,3)--(7,3);  
  \draw[thick,dotted](-1.5,4.5)--(4.5,4.5)--(4.5,1.5)--(-1.5,1.5)--cycle;  
  \draw[thick,dotted](5.5,4.5)--(11.5,4.5)--(11.5,1.5)--(5.5,1.5)--cycle;    
  \draw (.5,4.8) node {\emph{inversions }};
  \draw (10,4.8) node {\emph{knots}};
\end{tikzpicture}
\end{center}
\caption{} %
\label{fig:phylo}   
}
\end{figure}

If a unified picture of bacterial evolution is constructed using the algebraic ideas contained here, then one consequence will be that the tremendously rich theory behind finite reflection groups and their $q$-analogues, the Iwahori-Hecke algebras, as well as computational algebra systems such as GAP and Magma~\citep{GAP4,bosma1997magma}, will be made available to biologists as another powerful resource for their models.  Algebraic models can be expected to throw up new suggestions about the way biological mechanisms behave, and this can lead to new hypotheses for biologists to test.  Exactly this occurred in the case of a predicted knot through tangle algebra applied to the Tn3 recombinase~\citep{wasserman1985discovery}. 

A central problem of modern biology is the challenge of dealing with very large volumes of data.  For this reason, statistical methods have been the main plank of the biologists mathematical scaffolding.  An algebraic approach looks at the structures in a way that while not replacing statistical approaches, provides a new angle to tackling large volumes of information. By modelling using group theory and other algebraic structures, the extensive algebraic results and sophisticated and efficient algorithms of computational algebra can be brought to bear.  New computational and bioinformatic tools can then be developed to aid biologists.

From the algebraic side, the connections that are developed between algebraic structures and real biological questions will help to motivate further research in the algebraic structures themselves, and will raise new questions for algebraists.   Algebraists are familiar with structures such as reflection groups and Iwahori-Hecke algebras having a wide range of applications in certain parts of science.  For instance, reflection groups are natural ways to study symmetries arising in nature (e.g. crystallography), and it will surprise no algebraist to learn that they have already been used on occasion to count genetic arrangements using Burnside's Lemma (e.g.~\citet{ancelmeyers2005evolution}).   Similarly, Iwahori-Hecke algebras are important in the study of quantum groups, which first arose through quantum physics.  These algebraic structures arise widely because of their fundamental links to symmetries and patterns that arise in many places.  Perhaps they may find a greater role in biology as interdisciplinary work evolves.

\section{Acknowledgements}
I would like to thank Mark M. Tanaka, Leonard L. Scott Jr, John J. Graham and Attila Egri-Nagy, who read and commented on the manuscript.  Particular thanks to MMT who introduced me to the field of bacterial genomics in the first place.

\bibliographystyle{plainnat}

\end{document}